\documentclass[12pt]{article}
\usepackage[utf8]{inputenc}
\usepackage{amsmath}
\usepackage{xcolor}
\usepackage{tikz-cd}
\usepackage{amssymb}
\usepackage{amsthm}
\usepackage{esvect}
\usepackage{mathtools}
\usepackage{harpoon}

\usepackage[shortlabels,inline]{enumitem}
\usepackage{graphicx}
\usepackage[justification=centering]{caption}
\usepackage{subcaption}
\usepackage{comment}
\usepackage[T1]{fontenc}
\usepackage{float}
\usepackage{xcolor}
\newcommand{\FC}{F\mathcal C}
\newcommand{\states}{\mathcal V}
\newcommand{\cross}{\pitchfork}
\newcommand{\X}{\X}

\newcommand{\parl}{\mathbin{\|}}
\newcommand{\parw}{P_W^{\parl}}
\newcommand{\crossw}{P_W^\cross}
\newcommand{\stabw}{P_W^|}
\newcommand{\back}{B_W}
\newcommand{\forw}{A_W}
\renewcommand{\X}{\widetilde X}
\newcommand{\defeq}{\vcentcolon=}

\DeclareMathOperator{\Aut}{Aut}
\DeclareMathOperator{\Fix}{Stab}

\usepackage{graphicx}
\usepackage{comment}
\usepackage[T1]{fontenc}
\usepackage{float}
\usepackage{xcolor}

\newtheorem{theorem}{Theorem}[section]
\newtheorem{lemma}[theorem]{Lemma}
\newtheorem{corollary}[theorem]{Corollary}
\newtheorem*{conjecture}{Conjecture}
\newtheorem*{mtheorem}{Main Theorem}

\theoremstyle{definition}
\newtheorem{definition}[theorem]{Definition}

\usepackage{fullpage}

\title{\Large {Random Group Actions on $\mathrm{CAT}(0)$ Square Complexes}}
\author{\normalsize Zachary Munro}
\date{}

\begin{document}

%
%
%
%
%
%
%
%

\maketitle

\begin{abstract}
	\centering
	\begin{minipage}{\dimexpr\paperwidth-10cm}
		Generalizing ideas in \cite{jahncke}, we introduce the notion of progression in $\mathrm{CAT}(0)$ square complexes. Using progression, we are able to build on the proof strategy of \cite{dontsplit} to show any action of a random group with seven or more generators on a $\mathrm{CAT}(0)$ square complex has a global fixed point.
	\end{minipage}
\end{abstract}

\section{Introduction}
In this article, we begin the inquiry into the cubical dimension of random groups in the Gromov density model, first introduced in \cite{gromov}. Letting $S=\{s_1,\cdots ,s_n\}$ be a set of generators, a \emph{random group at density $d\in (0,1)$ and length $L$} is given by the presentation $\langle S\ |\ R\rangle$, where $R$ is a set of $\lfloor (2n-1)^{dL}\rfloor$-many words chosen uniformly at random from the set of reduced words of length $\ell$ with alphabet $S$. One says that a random group has some property $Q$ \emph{with overwhelming probability (w.o.p.)} if the probability that $\langle S\ |\ R\rangle$ satisfies $Q$ approaches 1 as $L\to \infty$. A number of results in the theory of random groups describe ``phase transitions'' for particular properties $Q$ with respect to the density $d$. It is proven in \cite{Ollivier}, \cite{gromov} that for $d>1/2$ a random group is either trivial or $\mathbb{Z}_2$ w.o.p., and for $d<1/2$ a random group is infinite, hyperbolic, and torsion-free w.o.p. Thus most of the interesting theory occurs in the range $d\in (0,1/2)$. 

There have been various results regarding phase transitions associated with actions of a random group on a $\mathrm{CAT}(0)$ cube complex. In the positive direction, in which one produces actions on $\mathrm{CAT}(0)$ cube complexes, Ollivier and Wise $\cite{ollivierwise}$ showed that random groups with density $d<\frac 15$ act without global fixed point on $\mathrm{CAT}(0)$ cube complexes of finite dimension w.o.p. This result has been strengthened by MacKay and Przytycki $\cite{balanced}$, who established the result for $d<\frac {5}{24}$, and Montee \cite{Montee}, who improved the bound to $d<\frac{3}{14}$. Recently, Ashcroft \cite{Ashcroft} has proven that random groups at densities $d<1/4$ act without global fixed point on $\mathrm{CAT}(0)$ cube complexes of finite dimension w.o.p. In the negative direction, it was shown by \.Zuk \cite{zuk} and Kotowski-Kotowski \cite{kotowski} that for $d>1/3$ a random group has Kazhdan's property (T) w.o.p. Groups with property (T) act always with global fixed point on $\mathrm{CAT}(0)$ cube complexes, following from an observation by Niblo and Reeves \cite{nibloreeves}. The following conjecture was communicated to the author by Piotr Przytycki \footnote{Note that $\frac 15$  is $\frac 45\cdot \frac 14$, $\frac {5}{24}$ is $\frac 56\cdot\frac 14$, and $\frac {3}{14}$ is $\frac 67\cdot\frac 14$.}. 

\begin{conjecture}
	Random groups act without global fixed point on finite dimensional $\mathrm{CAT}(0)$ cube complexes at densities $d<\frac 14$ and have property (T) at densities $d>\frac 14$ w.o.p.
\end{conjecture}

Note that Ashcroft's work in \cite{Ashcroft} completes the first half of the conjecture. 

For the purposes of this article, we concern ourselves with the dimension of those cube complexes on which random groups  act without global fixed point. In \cite{dontsplit}, it was shown that random groups at any density act always with global fixed point on 1-dimensional $\mathrm{CAT}(0)$ cube complexes. We extend this result by proving the following. 

\begin{mtheorem}
Random groups $G=\langle S\mid R\rangle$ with $|S|\geq 7$ act always with global fixed point on 2-dimensional $\mathrm{CAT}(0)$ cube complexes w.o.p.
\end{mtheorem}

A natural way of proving that a random group $\langle S\ |\  R\rangle$ does not exhibit a particular property $Q$ is to show there exists a large collection of words $\mathcal L$ in the alphabet $S$ which must be nontrivial if $\langle S\ |\ R\rangle$ were to have property $P$. Because $\mathcal L$ is large and $R$ consists of randomly sampled words, the probability that the intersection $\mathcal L\cap R$ is nonempty approaches 1 as $L\to \infty$. And thus we conclude that w.o.p. $\langle S\ |\ R\rangle$ does not have property $P$. This approach was applied fruitfully in each of \cite{dontsplit}, \cite{orlef}, and \cite{jahncke} to prove random groups do not have particular properties. We are especially inspired by \cite{jahncke}, where it is proved random groups act always with global fixed point on $\mathbb R$-trees.  

In Section 2 we establish our notation regarding $\mathrm{CAT}(0)$ cube complexes and graphs. In Section 3 we define checkpoint automata, a central tool in the paper. We prove a lemma ensuring that the accepted language of a checkpoint automaton is large. Given a group action $G\to \Aut(X)$ on some set $X$, an element $g\in G$ acts nontrivially if its image in $\Aut(X)$ is nontrivial. The lemmas proven in Section 4 are used to ensure that many elements of a group action are nontrivial. It is here that we prove the central technical theorem. We are immediately able to deduce a weakened form of our main theorem (Corollary~\ref{cor:weakenedMain}), where the number of generators depends on the density $d$. Section 5 is dedicated to removing the dependency on $d$, using a modified version of a construction found in \cite{dontsplit}. Given a random group $G=\langle S \mid R\rangle$ at some density $d$, one finds a finite index subgroup $H<G$ which is a quotient of a random group at density $d$ with exponentially more generators than $G$. Property $\FC_2$ is inherited by quotients by Lemma~\ref{lem:FCquotient} and $\FC_2$ is a commensurability invariant by Lemma~\ref{lem:FCvirtual}, and thus we deduce $G$ has property $\FC_2$.

\textbf{Acknowledgements.} 
I would like to express my immense gratitude to my advisor Piotr Przytycki, without whom none of this work would be possible. Thank you for posing this problem to me. Thank you for our regular meetings. And thank you for encouraging me when the math seemed hopeless. Thank you Adrien Abgrall for revising early drafts of this paper.

\section{Cube complexes and Graphs} 
We introduce our notation concerning $\mathrm{CAT}(0)$ square complexes and graphs below. For a more thorough introduction to $\mathrm{CAT}(0)$ cube complexes, see \cite{sageev}.

\subsection{Actions on cube complexes}

In our article, $\X$ will denote a $\mathrm{CAT}(0)$ square complex, and $S$ will be a set of formal letters for which there exists a homomorphism $F_S\to \Aut(\X)$, where $F_S$ is the free group on $S$. For a hyperplane $H$ of $\X$, we let $N(H)$ denote the carrier of $H$. We let $d(\cdot,\cdot)$ denote the combinatorial metric on $\X^1$. For 0-cells $x,y\in \X^0$, we let $[x,y]$ denote a combinatorial geodesic from $x$ to $y$. Note that $(\X^1,d)$ is not uniquely geodesic. Thus, if we ever assert that $[x,y]$ has some  property, then it is implicit that the property holds independent of the choice of geodesic. The intersection of a hyperplane with $\X^1$ is a set of midpoints of 1-cells in $\X^1$. The complement $\X-H$ of a hyperplane has two components, called \emph{halfspaces}. A hyperplane $H$ \emph{separates} subsets $A,B\subset \X$ if $A$ and $B$ are contained in distinct components of $\X-H$. If $H$ is a hyperplane and $x\in \X^0$, we let $d(x,H)\defeq \inf\{d(x,h)\ \mid \ h\in H\}$ denote the smallest distance from $x$ to $H$. For $x\in \X^0$ recall that $d(x,H)>1$ if and only if $x\not\in N(H)$ if and only if $x$ is separated from $H$ by some other hyperplane. Also, for $x,y\in \X^0$ the distance $d(x,y)$ is equal to the number of hyperplanes which separate $x$, $y$. Two distinct hyperplanes $H\neq W$ \emph{cross}, denoted by $H\cross W$, when $H\cap W\neq \emptyset$. If we write $H\cap W\neq\emptyset$, we leave open the possibility that $H=W$. Two distinct hyperplanes $H\neq W$ are \emph{parallel} if $H\cap W=\emptyset$. This is denoted $H\parl W$. 

An isometry $s$ of a $\mathrm{CAT}(0)$ cube complex is a \emph{hyperplane inversion} if it stabilizes some hyperplane $H$ and interchanges its halfspaces. In this case, we say $s$ \emph{inverts} $H$. Any group action is without hyperplane inversions after a subdivision of the cube complex. 


We will make use of 2-dimensionality in a couple of ways. For one, no more than two hyperplanes can pairwise intersect. Secondly, if $H\cross W$, then $N(H)\cap N(W)=C$ where $C$ is a single square of $\X$. In particular, for any 0-cell $x\in N(H)-C$ there is a hyperplane separating $x$ and $W$. 

A group $G$ has \emph{property $\FC_n$} if every action of $G$ by cubical isometries on an $n$-dimensional $\mathrm{CAT}(0)$ cube complex has a global fixed point. $G$ has property $\FC_\infty$ if $G$ has $\FC_n$ for every $n\in \mathbb N$. Note that $\FC_{n+1}\subset \FC_n$ and that $\FC_1$ is equivalent to Serre's property $\mathrm{FA}$. Similar to property $\mathrm{FA}$, the property $\FC_n$ passes to quotients since group actions can always be pulled back along group homomorphisms. Consequently, we have the following.

\begin{lemma}
\label{lem:FCquotient}
	If $G$ is a quotient of some $\widehat G$ with property $\FC_n$, then $G$ has property $\FC_n$.
\end{lemma}

The commensurability invariance of property $\mathrm{FA}$ also generalizes to $\FC_n$. A finite set of points in a $\mathrm{CAT}(0)$ cube complex has a unique barycentre, i.e. a point which minimizes the sum of distances to all points in the set. The existence of barycentres immediately leads to the following lemma. 

\begin{lemma}
\label{lem:FCvirtual}
	If a finite index subgroup $H<G$ has property $\FC_n$, then $G$ has property $\FC_n$. 	
\end{lemma}

Note the converse of the above lemma does not hold in general. 

\subsection{Graphs}

A \emph{graph} $\Gamma$ is a 1-complex whose 0-cells $V(\Gamma)$ are called \emph{vertices} and 1-cells $E(\Gamma)$ are called \emph{edges}. A map between complexes is \emph{combinatorial} if it sends $n$-balls homeomorphically to $n$-balls. An \emph{immersion of graphs} $\Gamma\looparrowright \Gamma'$ is a locally injective combinatorial map. If $P$ is a subdivided interval $[0,n]$ with vertices at integer points, then a combinatorial map $P\to \Gamma$ is a \emph{path of length $n$}. We think of paths as being oriented from 0 to $n$. The \emph{initial vertex} is the first vertex in the image of $P$ and the \emph{terminal vertex} is the last vertex in the image of~$P$. An \emph{initial segment} of $P$ is the restriction of $P\to \Gamma$ to some $[0,i]$ for $0\leq i \leq n$. 

Given a combinatorial map $\Gamma\to \Gamma'$ one can iteratively identify incident edges of $\Gamma$ with the same image. This procedure allows one to factor $\Gamma\to \Gamma'$ as a sequence $\Gamma_0\to \cdots \Gamma_n\to \Gamma'$ of combinatorial maps where $\Gamma_0=\Gamma$, and $\Gamma_i\to \Gamma_{i+1}$ is $\pi_1$-surjective for $0\leq i\leq n-1$, and $\Gamma_n\to \Gamma'$ is an immersion. The graph $\Gamma_n$ is uniquely determined as the coarsest quotient of $\Gamma$ such that the map $\Gamma\to\Gamma'$ factors as $\Gamma\to\Gamma_n\looparrowright\Gamma'$. We say that one obtains $\Gamma_n$ by \emph{folding}~$\Gamma$. 

A connected graph $T$ is a \emph{tree} if $\pi_1 T=\{1\}$. Note that a tree is uniquely geodesic. is A \emph{leaf} $v$ of a tree is a vertex which lies on the boundary of a single edge. A \emph{rooted tree} $(T,x)$ is a a tree with a distinguished vertex $x\in V(T)$. A vertex $v$ in a rooted tree $(T,x)$ has \emph{depth $m$} if the distance from $v$ to $x$ is $m$. The \emph{depth} of a rooted tree is the maximal depth across all vertices. In a rooted tree $(T,x)$, a vertex $u$ is an \emph{ancestor} of $v$ if $u$ lies on the geodesic from $x$ to $v$. The vertex $v$ is a \emph{descendant} of $u$ if $u$ is its ancestor. A \emph{child} of $u$ is a descendent which shares an edge with $u$. A rooted tree is a \emph{k-child} tree if each vertex has exactly $k$ children or is a leaf. 


\section{Large languages and checkpoint automata}

\subsection{Growth rate}

\begin{definition}
\label{def:wop}
	Let $S=\{s_1,\ldots, s_n\}$ be a set of formal letters and $F_S$	 the free group generated by $S$. For an integer $L>0$, let $R_L$ be the set of all reduced words of length $L$ over $S$. A \emph{random set of relators at density $d$ and length $L$} is a $\lfloor (2n-1)^{dL}\rfloor$-tuple elements from $R_L$ selected uniformly at random. A \emph{random group at density $d$ and length $L$} is defined by a presentation $\langle S\mid R\rangle$ where $R$ is a random set of relators at density $d$ and length $L$. Let $P$ be a property for groups (or a set of relators), let $d\in (0,1)$ be a density, and let $\mathbb{P}_L(P)$ be the probability that $P$ holds for a random group (or set of relators) at density $d$ and length $L$. Property $P$ holds \emph{with overwhelming probability} (w.o.p.) if $\mathbb P_L(P)\to 1$ as $L\to\infty$. 
\end{definition}

The following definition quantitatively captures what it means for a language to be large.

\begin{definition}
    Let $S=\{s_1,\ldots, s_n\}$ be a set of letters and $\mathcal L$ a language of reduced words over $S$. The language $\mathcal L$ is said to have \emph{growth rate at least k} if $\mathcal L\cap R_L>ck^L$ for some $c>0$ and all but finitely many $L>0$. 
    
\end{definition}

The following lemma first appears in Gromov's article \cite{gromov}. A detailed proof can be found in \cite{orlef}. 

\begin{lemma}
\label{lem:intersection}
	A language $\mathcal L$ of reduced words with growth rate $k>(2n-1)^{1-d}$ intersects a random set of relators $R$ at density $d$ w.o.p. 	
\end{lemma}

The following corollary about the intersection of a large language and a set of random relators was proven in \cite{dontsplit}.

\begin{corollary}[\cite{dontsplit}, Lemma 2.4]
\label{cor:intersection}
	Suppose $\lceil \lambda 2n\rceil -1 > 0$ and $(2n-1)^d\geq \frac2\lambda$ for some $\lambda\in (0,1)$. If a language $\mathcal L$ of reduced words has growth rate at least $(\lceil \lambda 2n\rceil -1)$, then a random set of relators $R$ at density $d$ intersects $\mathcal L$ w.o.p.
\end{corollary}

Let $S$ and $F_S$ be as in Definition~\ref{def:wop}. Suppose we are given an action $F_S\to \Aut(\X)$ on a $\mathrm{CAT}(0)$ cube complex. Then reduced words $w\in S^*$ denote an element of $F_S$, can be identified with an automorphism of $\X$ via the map $F_S\to \Aut(\X)$. The key technical theorem in this article is the following. 

\begin{theorem}
\label{thm:main1}
		There exists a finite set of languages $\mathcal L_{\Sigma^1}, \ldots, \mathcal L_{\Sigma^k}$ with $\frac n6$-growth such that the following holds: For any $\mathrm{CAT}(0)$ square complex and any action $F_S\to \Aut(\X)$ without global fixed point, there exists $j\in\{1,\ldots,k\}$ so that every $w\in \mathcal L_{\Sigma^j}$ acts nontrivially on $\X$.
\end{theorem}

Combining Corollary~\ref{cor:intersection} with the above theorem we can immediately deduce $\FC_2$ for certain random groups. Taking $\lambda=\frac{1}{12}$, we get that at each density $d$, random groups with sufficiently large generating set have property $\FC_2$.

\begin{corollary}
\label{cor:weakenedMain}
	Fix $0<d<1$ and $n\geq 7$. Let $S=\{s_1,\ldots, s_n\}$ be a set of formal letters such that $(2n-1)^d\geq 24$. Then a random group $\langle S \mid R\rangle$ at density $d$ does not act without global fixed point on a $\mathrm{CAT}(0)$ square complex w.o.p.
\end{corollary}

\subsection{Checkpoint automata}

Let $S=\{s_1,\ldots s_n\}$ freely generate $F_S$. Let $B_S$ be a bouquet of directed circles labeled with elements of $S$ so that $\pi_1 B_S=F_S$. If the edges of a graph $\Gamma$ are directed and labeled by elements of $S$, then there is a combinatorial map $\Gamma\to B_S$ induced by the directions and labels. 

\begin{definition}
\label{def:languagetree}
    A \emph{checkpoint automaton (c-automaton)} $\Sigma$ is a directed graph with an immersion $\Sigma\looparrowright B_S$ constructed as follows.
    
    \begin{enumerate}
    	\item Begin with a discrete set of \emph{checkpoint vertices} $\states$ with a distinguished \emph{start vertex} $v_0\in \states$.
    	\item To each vertex $u\in \states$ we associate a directed, $S$-labeled, finite rooted tree $T_u$ with root $u$ and $E(T_u)\neq\emptyset$.
    	\item Let $l(T_u)$ denote the set of leaves of $T_u$. For each $u\in \states$, we label $l(T_u)$ with elements of $\states'\defeq \states-\{v_0\}$.
    \end{enumerate}

The c-automaton $\Sigma$ is the quotient of $\sqcup_{u\in \states} T_u$ induced by the labels on $\{l(T_u)\}_{u\in \mathcal V}$. See Figure~\ref{fig:AutomatonData}. We require that the natural map $\Sigma\to B_S$ is an immersion. The \emph{accepted language} of $\Sigma$ is the set of words $\mathcal L_\Sigma$ which label immersed paths $P\looparrowright \Sigma$ beginning at $v_0$ such that for each $u\in \states$ any subpath in $P\cap T_u$ is directed away from the root $u\in T_u$ . See Figure~\ref{fig:cAutomaton}.

An \emph{automaton} $\Sigma$ is a c-automaton with $\mathcal V= V(\Sigma)$ and $T_u$ has depth one for each $u\in \mathcal V$. Note that any c-automaton can be easily converted to an automaton by enlarging the set of checkpoint vertices. This operation does not affect the accepted language. Thus, there is a natural equivalence automata and c-automata. The structure of c-automata is merely a formalism enabling us to more easily express certain ideas.  

\begin{figure}[h]
	\centering
	\begin{subfigure}[h]{0.45\textwidth}
		\centering
        \includegraphics[width=\textwidth]{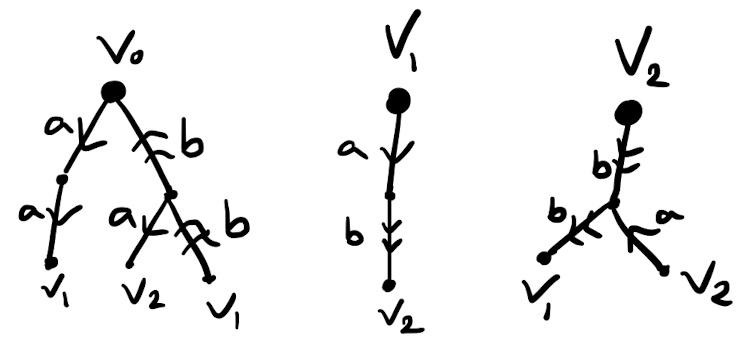}
        \caption{The data encoding a c-automaton with $\mathcal V=\{v_0,v_1,v_2\}$.}	
        \label{fig:AutomatonData}
	\end{subfigure}
	\ \ \ \ \ 
	\begin{subfigure}[h]{0.4\textwidth}
		\centering	
		\includegraphics[width=\textwidth]{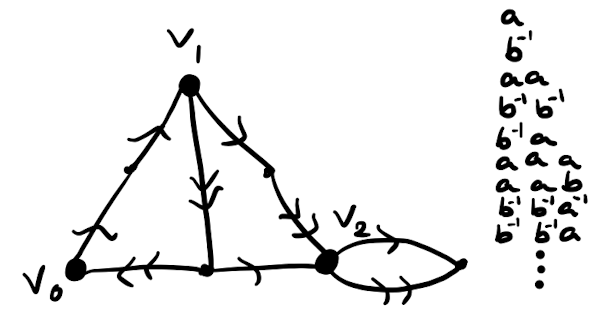}
		\caption{The c-automaton with some words in the accepted language $\mathcal L_\Sigma$.}	
		\label{fig:cAutomaton}
	\end{subfigure}
	\caption{}
\end{figure}

\end{definition}

We will make use of a local largeness criterion, depending only on the trees $T_u$, to guarantee that the accepted language $\mathcal L_\Sigma$ is large. 
\begin{definition}


A vertex $u\in \states$ has \emph{$k$-growth} if either $u=v_0$ and $T_{v_0}$ has depth one, or for every $v\in V(T_u)$ either $v$ is a leaf or has at least $k$ children. 
\end{definition}

The following lemmas follow from definitions. 

\begin{lemma}
\label{growth}
	If all vertices in a c-automaton $\Sigma$ have $k$-growth, then $\mathcal L_\Sigma$ has growth rate at least $k$.	
\end{lemma}

\begin{lemma}
\label{lem:reduceStrongGrowth}
	Suppose $\Sigma$ is a c-automaton whose vertices all have $k$-growth. For each $u\in \states$, remove up to $k'<k$ children from each $v\in V(T_u)$. Let $\Sigma'_S$ be the resulting c-automaton. Then each vertex of $\Sigma'_S$ has $(k-k')$-growth.
\end{lemma}

\section{Progressing in Square complexes}
\label{sec:progressing}

Let $S$ be a set of formal letters, $F_S$ the free group over $S$, and $F_S\to \Aut(\X)$ a homomorphism to the automorphism group of a $\mathrm{CAT}(0)$ square complex. We define an automaton $\Sigma$ over $S$ as follows.

\begin{definition}
\label{def:states}
	Let $x$ be a 0-cell of $\X$. For each $s\in S$, let $\mathcal H_s$ be the set of hyperplanes $H$ such that $H \cap[x,sx]\neq \emptyset$ and $x\in N(H)$. Equivalently, those $H$ intersecting $[x,sx]$ such that there is no hyperplane separating $x$ and $H$. Elements of $\mathcal H_s$ pairwise cross, and thus there can be at most two hyperplanes in each $\mathcal H_s$.
\end{definition}

The set of such hyperplanes $\mathcal H = \bigcup_{s\in S}\mathcal H_s$ along with a special start vertex $v_0$ will be the checkpoint vertices $\mathcal V=\mathcal H \cup \{v_0\}$ of our automaton $\Sigma_{S}$. 

\begin{definition}
Let $x$, $\mathcal H_s$, and $\mathcal V$ be as above. Any path $P=e_1\cdots e_n\looparrowright B_S$ immersing into the bouquet of circles $B_S$ corresponds to a reduced word $\widehat P=\hat e_1\cdots \hat e_n\in S^*$. Suppose $P$ joins vertices $H$ and $H'$ in $\mathcal V$. We say $P$ is \emph{progressing with respect to $x$} if the following hold


	\begin{enumerate}
		\item If $H=v_0$, then $\widehat P H'$ separates $x$ and $\hat e_1\cdots \hat e_i x$ for each $1\leq i\leq n$.
		\item If $H\neq v_0$, then we have that
		
			\begin{enumerate}[label*=\arabic*]
				\item $H$ does not separate $x$ and $\hat e_1\cdots\hat e_ix$ for each $1\leq i\leq n$, and
				\item Either $\widehat P H'\parl H$ and $\widehat P H'$ separates $x$ and $\widehat P x$, or $\widehat P H'=H$.
			\end{enumerate}
		
	\end{enumerate}
	
\end{definition}

\begin{figure}[h]
	\centering
	\begin{subfigure}[h]{0.45\textwidth}
		\centering
        \includegraphics[width=\textwidth]{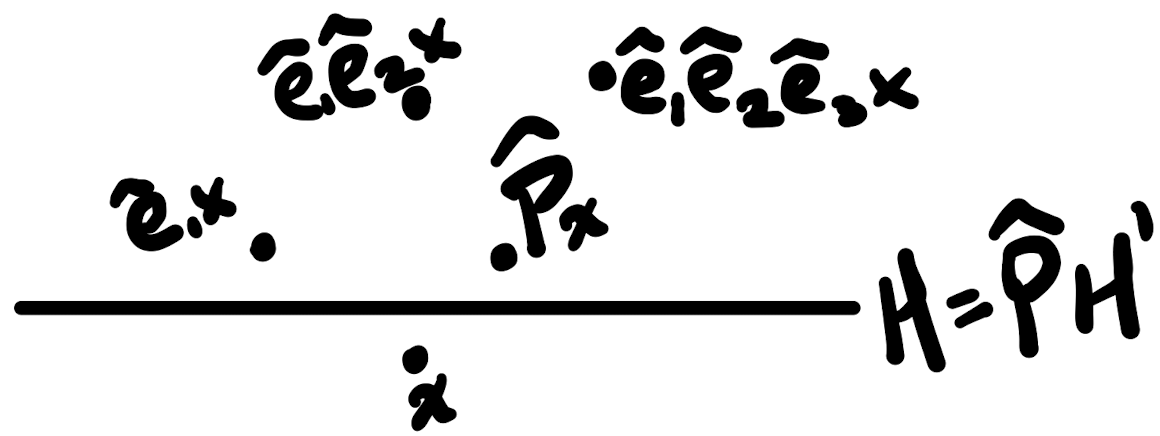}
        \caption{An example of a progressing $P$ satisfying $\widehat PH'=H$.}	
	\end{subfigure}
	\ \ \ \ \ 
	\begin{subfigure}[h]{0.4\textwidth}
		\centering	
		\includegraphics[width=\textwidth]{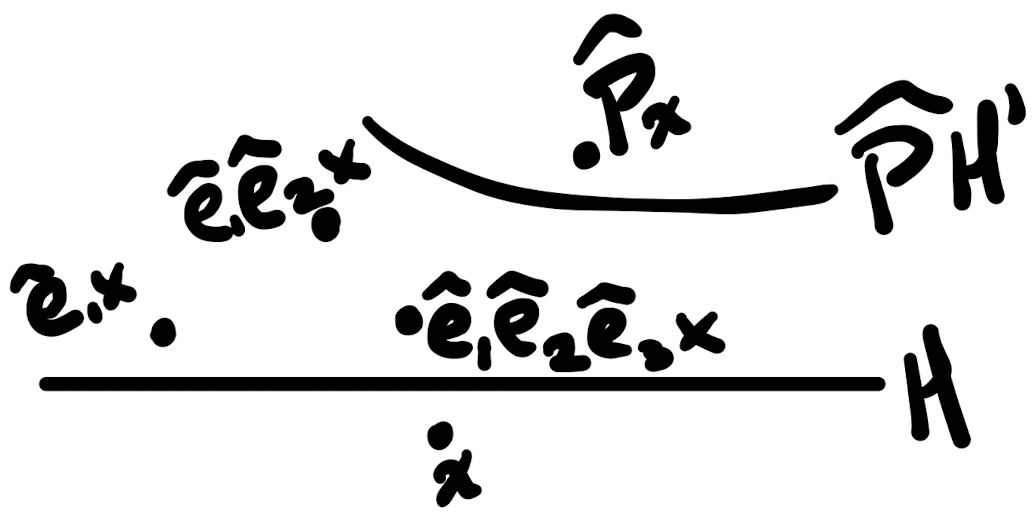}
		\caption{An example of a progressing $P$ satisfying $\widehat P H'\parl H$.}	
	\end{subfigure}
	\caption{}
\end{figure}

For $W\in \mathcal H$, a rooted tree $T_W$ as in Definition~\ref{def:languagetree}.2 is \emph{progressing with respect to $x$} if each path $P\looparrowright T_W$ joining $W$ to a leaf is progressing. A c-automaton $\Sigma$ is \emph{progressing with respect to $x$} if for $T_u$ is progressing with respect to $x$ for each $u\in \states$. By the definition of c-automata, there is a natural map $P\looparrowright T_u\looparrowright B_S$.

\begin{lemma}
If $\Sigma$ is progressing with respect to $x$, then $wx\neq x$ for each nontrivial accepted word $w\in \mathcal L_\Sigma$. In particular, the nontrivial accepted words act nontrivially on $\X$.   
\end{lemma}

\begin{proof}
We first prove the lemma for words labeling paths that end on checkpoint vertices. Suppose $w\in \mathcal L_\Sigma$ labels a path $P\looparrowright \Sigma$ beginning at $v_0$. Let $P=P^1\cdots P^n$ be the unique decomposition of $P$ into minimal subpaths joining two checkpoint vertices. Note minimality is equivalent to the existence of some $T_u$ for each $P^j$ such that $P^j\looparrowright T_u$ and joins the root to a leaf. Suppose these checkpoint vertices, in order, are $H_0,\ldots, H_n$. Let $e^i_1, \ldots, e^i_{n_i}$ be the edges which $P^i$ traverses in order. We prove by induction on $n$ that $\widehat PH_n$ separates $x$ and $\widehat P x$. 

Suppose $n=1$. Then $H_0=v_0$ and by definition $\widehat P H_1$ separates $x$ and $\widehat P x$.  

Let $P'\defeq P^1\cdots P^{n-1}$ and suppose $\widehat P' H_{n-1}$ separates $x$ and $\widehat P' x$. Since $\widehat P=\widehat P'\widehat P^n$ and $P^n$ is progressing, $\widehat P x$ lies on the same side of $\widehat PH_{n-1}$ as $\widehat P' x$. Thus $\widehat P H_{n-1}$ separates $x$ and $\widehat P x$. Moreover, $\widehat P H_n$ is either parallel to $\widehat P H_{n-1}$ and separates $\widehat P x$ and  $\widehat P' x$, or $\widehat P H_n=\widehat P'H_{n-1}$. In either case, we have that $\widehat P H_n$ separates $x$ and $\widehat P x$. 

The conclusion of the lemma is an immediate consequence of the above. Since $\widehat P^1\cdots \widehat P^i H_i$ separates $x$ and $\widehat P^1\cdots \widehat P^i x$, then by definition of progressing $\widehat P^1\cdots \widehat P^i H_i$ also separates $x$ and $\widehat P^1\cdots \widehat P^i\hat e^{i+1}_1\cdots \hat e^{i+1}_j x$ for all $1\leq j\leq n_{i+1}$ and $i>1$. If $i=1$, then the definition of progressing ensures $\hat e^1_1\cdots \hat e^1_j x$ and $x$ are separated by $\widehat P^1 H_1$ for all $1\leq j\leq n_1$. 

Observe we have proven that for any initial subpath $P'\subset P$, the points $\widehat P'x$ and $x$ are separated by a hyperplane. Since any accepted word can be extended to an accepted word that labels a path ending on a checkpoint vertex, we are done.

\end{proof}

\begin{figure}[h]
	\centering
	\begin{subfigure}[h]{0.3\textwidth}
		\centering
        \includegraphics[width=\textwidth]{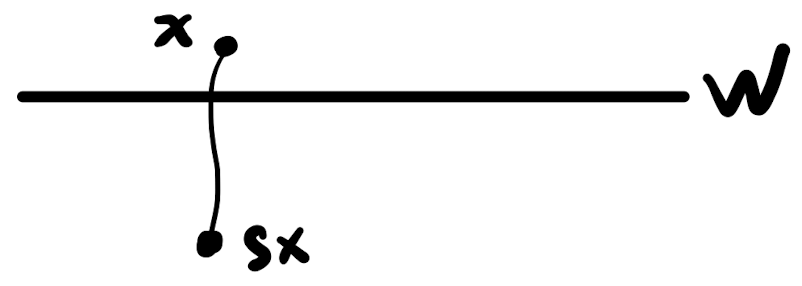}
        \caption{$s\in \back$}	
	\end{subfigure}
	\begin{subfigure}[h]{0.3\textwidth}
		\centering	
		\includegraphics[width=\textwidth]{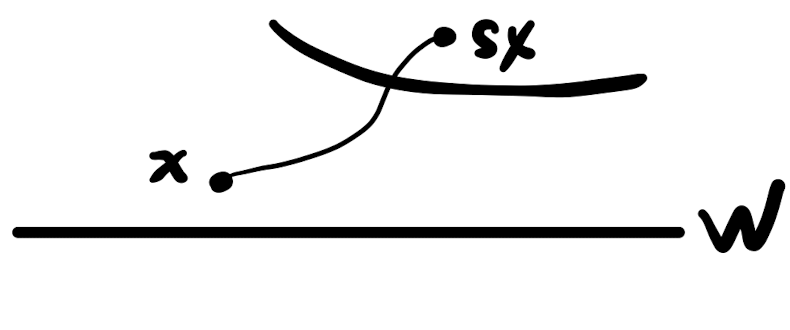}
		\caption{$s\in \forw$}	
	\end{subfigure}
	\begin{subfigure}[h]{0.3\textwidth}
		\centering
        \includegraphics[width=\textwidth]{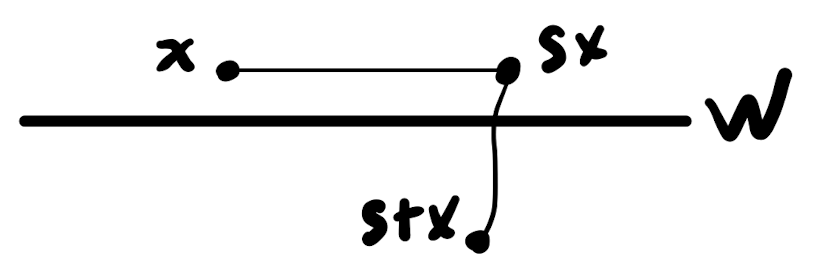}
        \caption{$s\in \stabw$}	
	\end{subfigure}
	
	\begin{subfigure}[h]{0.3\textwidth}
		\centering
        \includegraphics[width=\textwidth]{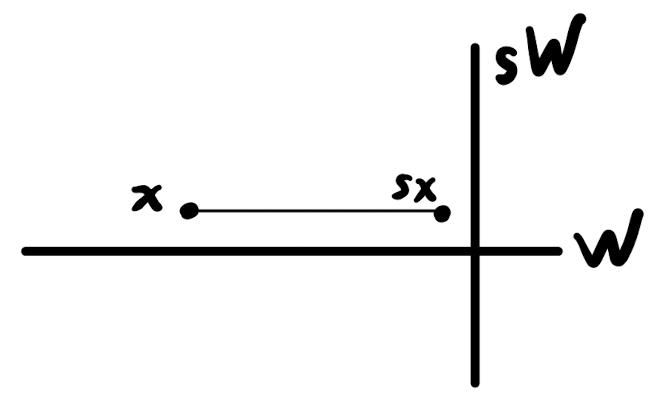}
        \caption{$s\in \crossw$}	
	\end{subfigure}
	\begin{subfigure}[h]{0.4\textwidth}
		\centering
        \includegraphics[width=\textwidth]{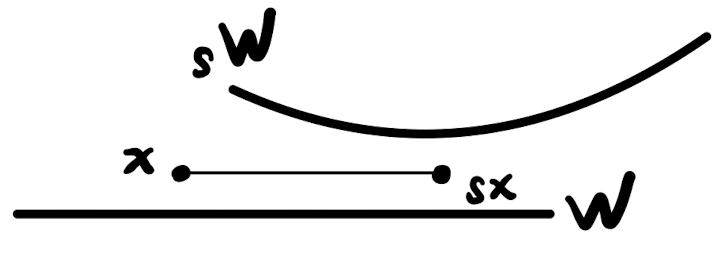}
        \caption{$s\in \parw$}	
	\end{subfigure}
	
	\caption{The partition from Definition~\ref{def:partition}}

	\label{fig:partition}
\end{figure}

\begin{definition}
\label{def:partition}
Fix $x\in \X^0$ and let $\mathcal H$ be as above. For each $W \in \mathcal H$ we define a partition of the letters and their inverses $S \cup S^{-1} = \forw \sqcup \back \sqcup \stabw \sqcup \crossw \sqcup \parw$ as follows. See Figure~\ref{fig:partition}.
	
	\begin{enumerate}
		\item $s\in \back$ if $[x,sx]$ crosses $W$
		
		\item $s\in \forw$ if $[x,sx]\cap W = \emptyset$ and $sx\not\in N(W)$
		
		\item $s\in \stabw$ if $[x,sx]\subset N(W)-W$ and $[sx,stx]\cap W\neq\emptyset$ for some $t\in S\cup S^{-1}$ 
		
		\item $s\in \crossw$ if $[x,sx]\subset N(W)-W$ and $[sx,stx]\cap W=\emptyset$ for all $t\in S\cup S^{-1}$ and $sW\cross W$
		
		\item $s\in \parw$ if $[x,sx]\subset N(W)-W$ and $[sx,stx]\cap W=\emptyset$ for all $t\in S\cup S^{-1}$ and $sW\parl W$
	\end{enumerate} 
\end{definition}

We set $P_W\defeq \stabw\cup \crossw\cup \parw$, which can be characterized as those $s\in S\cup S^{-1}$ such that $[x,sx]\subset N(W)-W$. Note that the stabilizer of $x$ is partitioned across the various sets above. In particular, $s\in \Fix (x)$ can belong to any of $\stabw$, $\crossw$, or $\parw$. Consequently, $P_W-\Fix (x)$ are exactly the $s\in S\cup S^{-1}$ such that $sx\neq x$ and $[x,sx]\subset N(W)-W$.

\begin{lemma}
\label{lem:progressing}
    Consider a path $P\looparrowright B_S$ with initial vertex $W\neq v_0$ and traversing the edges $e_1,\ldots, e_n$ in order. We aim to choose a terminal vertex $H$ which makes $P$ progressing.
    
    \begin{enumerate}
        \item If $n=1$ and $\hat e_1\in A_W$, then there exists a choice of $H\in \mathcal H_{\hat e_1^{-1}}$ which makes $P$ progressing.
        
        \item If $n=1$ and $\hat e_1\in \stabw$ and $t\in S\cup S^{-1}$ is such that $[\hat e_1 x, \hat e_1 t x]\cap W\neq \emptyset$, then there is a choice of $H\in \mathcal H_t$ which makes $P$ progressing.
        
        \item If $n=2$ and $\hat e_1\in \crossw$ and $\hat e_2\in P_W-\Fix (x)$, then there is a choice of $H\in \mathcal H_{\hat e_2^{-1}}$ which makes $P$ progressing.
        
        \item If $n=2$ and $\hat e_1\in \parw$ and $\hat e_2\in \back$, then there is a choice $H\in \mathcal H_{\hat e_2^{-1}}$ which makes $P$ progressing.
    \end{enumerate}
\end{lemma}

\begin{proof}
\ 
\begin{enumerate}
    \item By definition of $\forw$, there exists a hyperplane $H'\parl W$ separating $x$ and $\hat e_1 x$. If we pick $H'$ to be a closest such hyperplane to $\hat e_1 x$, then $\hat e_1^{-1} H'\in \mathcal H_{\hat e_1^{-1}}$. Setting $H\defeq \hat e_1^{-1} H'$ makes $P$ progressing at $x$.
   
    \item Since $\hat e_1\in \stabw$ the 0-cell $\hat e_1 x$ lies in $N(W)$. Then since $W$ separates $\hat e_1 x$ and $\hat e_1 t x$, the hyperplane $\hat e_1^{-1} W$ lies in $\mathcal H_t$. Setting $H\defeq \hat e_1^{-1} W$ makes $P$ progressing at $x$. 
    
    \item By definition of $\crossw$, we have  $\hat e_1x\in N(W)\cap N(\hat e_1 W)$, and $[\hat e_1 x, \hat e_1 t x]$ does not cross $W$ for any choice of $t\in S\cup S^{-1}$. For $\hat e_2\in (\stabw\cup \crossw\cup \parw)-\Fix (x)$, any geodesic $[x, \hat e_2 x]\subset N(W)$ has length greater than zero. Thus $[\hat e_1x, \hat e_1\hat e_2x]\subset N(\hat e_1 W)$ is cut by some hyperplane $H'\parl W$. Here we are using the fact that $N(W)\cap N(\hat e_1 W)$ is a single square, since $\X$ is 2-dimensional. Picking such an $H'$ closest to $\hat e_1\hat e_2 x$ we have $(\hat e_1\hat e_2)^{-1}H'\in\mathcal H_{\hat e_2^{-1}}$. Setting $H\defeq (\hat e_1\hat e_2)^{-1}H'$ makes $P$ progressing at $x$. 
    
    \item Since $[x, \hat e_2]$ crosses $W$ and $\hat e_1 W\parl W$, the hyperplane $\hat e_1W$ separates $\hat e_1 x$ and $\hat e_1\hat e_2 x$. If we pick a closest hyperplane $H'\parl W$ to $\hat e_1 \hat e_2 x$ separating it from $\hat e_1 x$, then $(\hat e_1\hat e_2)^{-1}H'\in \mathcal H_{\hat e_2^{-1}}$. Setting $H\defeq (\hat e_1\hat e_2)^{-1}H'$ makes $P$ progressing at $x$.
\end{enumerate}
\end{proof}

\begin{corollary}
\label{cor:progressing}
	Let $W\in \mathcal H$ be a checkpoint vertex. Let $\{P_i\to B_S\}_{i\in I}$ be a collection of labeled paths of forms 1, 2, 3, or 4 in Lemma~\ref{lem:progressing}. Each $P_i$ has initial vertex $W$, so we can identify all paths at their initial vertex $T'=\sqcup_i P_i/\sim$. Fold $T'\to B_S$ to get $T_W\to B_S$. Then $T_W$ is progressing. 
\end{corollary}

As noted before, our strategy to prove Theorem~\ref{thm:main1} is to generate a set of c-automata with large growth whose accepted languages consist of words acting nontrivially on $\X$. The notion of progression will allow us to ensure the accepted languages act nontrivially. The difficulty now becomes to find trees $T_W$ such that each maximal $P\looparrowright T_W$ is progressing and such that $T_W$ has enough leaves to guarantee large growth. Lemma~\ref{lem:progressing} demonstrates there are many options for $P\looparrowright T_W$ when any of $A_W$, $\crossw$, or $\stabw$ are large. Carefully choosing the basepoint will allow us to have more control when we are not in one of these cases.

\begin{lemma}
\label{lem:half}
Fix $x\in\X^0$ minimizing $\Sigma_{s\in S\cup S^{-1}} d(x,sx)$, and partition $S\cup S^{-1}$ as in Definition~\ref{def:partition}. Let $n=|S|$. If $S$ acts on $\X$ without hyperplane inversions, then $|\back|\leq n$ for any $W\in \mathcal H$. 
\end{lemma}
\begin{proof}
Suppose there exists $W\in \mathcal H$ such that $|\back|> |S\cup S^{-1}|/2$. Let $y\in \X^0$ be the other endpoint of the edge $e$ dual to $W$ at $x$. For each isometry and its inverse $s,s^{-1}\in S\cup S^{-1}$, we have the following three possibilities: 

\begin{enumerate}
    \item Both $s,s^{-1}\in \back$. We show that $d(y,sy)=d(x,sx)-2$. Let $x_0,\ldots, x_n$ be the 0-cells of a geodesic $P\to\X$ joining $s^{-1}x=x_0$ to $x=x_n$. Since $W$ separates $s^{-1}x$ and $x$, we can take $x_{n-1}=y$. If $d(s^{-1}x,x)>2$, then we can pick $x_1$ so that $sx_1=y$. Thus $d(y,sy)=d(s,sx)-2$. If $d(s^{-1}x,x)=1$, then $W$ is dual to the edge with endpoints $s^{-1}x,x$ and thus $s$ inverts $W$, a contradiction. If $d(s^{-1}x,x)=2$, then $x$, $sy$, $sx$ are the vertices of a geodesic. If $sy\neq y$, then $s$ inverts $W$, a contradiction. If $sy=y$, then $d(sy,y)=d(sx,x)-2$.

    \item Exactly one of $s\in \back$ or $s^{-1}\in\back$. Suppose $s\in\back$. Consider the 0-cells $x_0,\ldots, x_n$ of a geodesic $P\to\X$ joining $x=x_0$ to $sx=x_n$ with $x_1=y$. Then $s^{-1}P$ is a geodesic connecting $s^{-1}x$ to $x$. The hyperplane $W$ does not intersect $s^{-1}P$, otherwise $W$ would separate $s^{-1}x$ and $x$. Thus the path formed by removing the last edge of the concatenation $e\cdot s^{-1}P$ is a geodesic joining $y$ to $s^{-1}y$. Since this geodesic is the same length as $P$, we have $d(y,sy)=d(x,sx)$. This equality is also true when $s^{-1}x\in \back$ by a symmetric argument.

    \item Both $s,s^{-1}\not\in \back$. Since $d(x,y)=1$, we have that $d(y,sy)\leq d(x,sx)+2$.
 
\end{enumerate}

Let the partition $S\cup S^{-1}=D_1\sqcup D_2\sqcup D_3$ correspond to the three cases above. Since $|\back|> |S\cup S^{-1}|/2$, we have that $|D_1|>|D_3|$. This yields
\[
\sum_{s\in S\cup S^{-1}}d(x,sx)\geq\sum_{s\in D_1}\big( d(y,sy)+2 \big) +\sum_{s\in D_2}d(y,sy) + \sum_{s\in D_3}\big(d(y,sy)-2\big)>\sum_{s\in S\cup S^{-1}}d(s,sy),
\]
contradicting the minimality of $\Sigma_{s\in S\cup S^{-1}} d(x,sx)$. 
\end{proof}

We our now ready to prove our main theorem. Let $n\defeq |S|$ so that $2n=|S\cup S^{-1}|$.

\begin{proof}[Proof of Theorem~\ref{thm:main1}]
	Let $x\in \X^0$ be a point minimizing $\Sigma_{s\in S\cup S^{-1}} d(x,sx)$, and let $\mathcal H$ be as in Definition~\ref{def:states}.
		
	\textbf{Case 1:} $|\Fix(x)|\geq \frac n6$
	
	If $|\Fix (x)|\geq \frac n6$, then it is possible to construct a c-automaton $\Sigma$ with $\frac n6$-growth. We set $\states=\{v_0\}$ and add a directed loop at $v_0$ with label $s$ for each $s\in \Fix(x)$. It is assumed in the remaining cases that $|\Fix(x)|<\frac n6$. 
	
	\textbf{Case 2:} For each $W\in \mathcal H$ one of the following hold:
	
	\begin{enumerate}
		\item $|\forw|\geq\frac n6$
		\item $|\stabw| \geq\frac {5n}{18}$
		\item $|\crossw|\geq \frac {5n}{18}$
		\item $|\back| \geq \frac n6$ and $|\parw|\geq \frac{5n}{18}$
	\end{enumerate}
	
	For each of the cases 1, 2, 3, and 4, we construct a progressing tree $T_W$ with root $W$ such that $T_W$ has $\frac n6$-growth. Since each $W$ satisfies one of 1, 2, 3, and 4, our various $T_W$ for $W\in \mathcal H$ form a c-automaton with $\frac n6$-growth. 
	
	Suppose that $|\forw| \geq \frac n6$. Then by Lemma~\ref{lem:progressing}.1, it is possible to construct $T_W$ with $\frac n6$-growth. 
	
	If $|\stabw|>\frac{5n}{18}>\frac n6$, then by Lemma~\ref{lem:progressing}.2 we can construct a progressing $T_W$ with $\frac n6$-growth.  
	
	We have $|P_W| > \frac{5n}{6}$ because $|\back| \leq n$ by Lemma~\ref{lem:half}. If $|\crossw|>\frac{5n}{18}$, then since $|P_W-\Fix (x)| > \frac{5n}{18}-\frac n6 = \frac {2n}{3}$ we can construct a progressing $T_W$ with $\frac n6$-growth by Lemma~\ref{lem:progressing}.3. 
	
	By Corollary~\ref{cor:progressing}, $T_W$ can be constructed to have $|\parw|\cdot|\back|$-many leaves of depth two. Thus we achieve $\frac n6$-growth if both $|\parw| >\frac {5n}{18}$ and $|\back| \geq \frac n6$.

	\textbf{Case 3:} There exists a $W\in \mathcal H$ such that all of the following hold:
	
	\begin{enumerate}
		\item $|\forw|<\frac n6$
		\item $|\stabw| <\frac {5n}{18}$
		\item $|\crossw|<\frac {5n}{18}$
		\item $|\back| < \frac n6$ and $|\parw|>\frac {5n}{18}$
	\end{enumerate}

	Suppose that $|\forw| < \frac n6$ and $|\Fix (x)| < \frac n6$. Thus one of $|\stabw| >\frac {5n}{18}$, $|\crossw| > \frac {5n}{18}$, or $|\parw| > \frac {5n}{18}$ must hold. Thus the above three cases are exhaustive. 
	
	Since $|\forw| < \frac n6$, we have $|P_W-\Fix(x)|>\frac{3n}{2}$. The $P_W$ translates of $x$ lie in the same component of $N(W)-W$. This implies that $|\mathcal H_s|=1$ for each $s\in P_W-\Fix(x)$. Indeed, any $H\in \mathcal H_s$ crosses $W$ and hyperplanes in $\mathcal H_s$ pairwise cross. No three hyperplanes can pairwise cross since $\dim (\X)=2$ and thus $|\mathcal H_s|=1$. We hereby identify $\mathcal H_s$ with the single hyperplane it contains for each $s\in P_W-\Fix(x)$. 
	
	If $\mathcal H_s\neq\mathcal H_t$ for $s,t\in P_W-\Fix(x)$, then $\mathcal H_s \parallel \mathcal H_t$ since both hyperplanes cross $W$. This implies $t\in A_{\mathcal H_s}$ and $s\in A_{\mathcal H_t}$. See Figure~\ref{fig:flowerPetals}. Fixing $s\in P_W-\Fix(x)$, we have $|A_{\mathcal H_s}\cap P_W|=|(P_W-\Fix(x)) - B_{\mathcal H_s}|\geq\frac{3n}{2}-n=\frac n2$. The set of generators $\mathcal A_{\mathcal H_s}\defeq \{t\in A_{\mathcal H_s}\cap P_W :\ t^{-1}\in (P_W-\Fix(x)) \}$ contains at least $\frac n6$-many elements since $|A_W|,|B_W|\leq \frac n6$. 
	
	\begin{figure}[h]
		\centering
		\includegraphics[width=6cm]{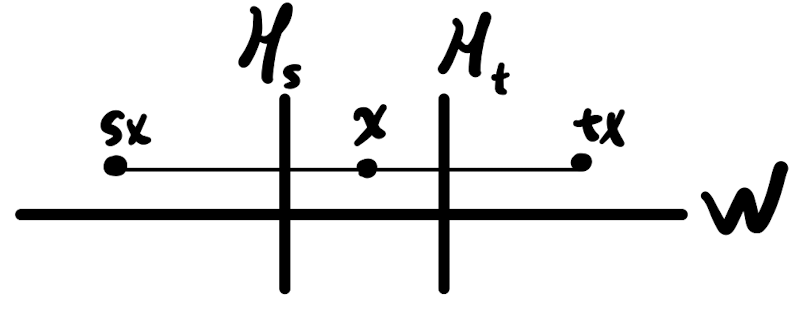}
		\caption{$t\in A_{\mathcal H_s}$ and $s\in A_{\mathcal H_t}$}
		\label{fig:flowerPetals}
	\end{figure}
	
	We construct a progressing c-automaton $\Sigma$ with $\frac n6$-growth in this final case. We let $V(\Sigma)=\{v_0\}\cup \{\mathcal H_s\ :\ s\in P_W-\Fix(x)\}$. For each $s\in P_W-\Fix(x)$ we add edges labeled by $\mathcal A_{\mathcal H_s}$ of the form in Lemma~\ref{lem:progressing}.1. Picking an arbitrary $s\in P_W-\Fix(x)$, we add an edge with initial vertex $v_0$ and terminal vertex $\mathcal H_{s^{-1}}$. Since $|\mathcal A_{\mathcal H_s}|\geq \frac n6$ for each $s\in P_W-\Fix(x)$, the c-automaton $\Sigma$ has $\frac n6$-growth. By Lemma~\ref{lem:progressing}.1, $\Sigma$ is progressing.
	
\end{proof}

Note that $|\mathcal H_s|\leq 2$ for each $s\in S\cup S^{-1}$, so that $|\mathcal H|\leq 4n$. The construction in the proof of Theorem~\ref{thm:main1} thus produces a finite set of automata since there are at most $4n+1$ vertices in  $\states= \mathcal H\cup \{v_0\}$ and each tree $T_u$ for $u\in \states$ has depth at most two.

\section{Large generating sets and pruning}
\label{sec:generators}

Let $\widetilde S$ denote the set of reduced words of length $B>0$ with letters in $S\cup S^{-1}$. The involution on $S\cup S^{-1}$ defined by $s\mapsto s^{-1}$ extends to an involution on $\widetilde S$. Since this involution has no fixed points, it defines a partition $\widetilde S=\widehat S\cup \widehat S^{-1}$. We let $\hat n\defeq n(2n-1)^{B-1}$ denote the cardinality of $\widehat S$. Let $L$ denote the length of a random set of relators in the density model. Let $I_P$ be the set of $L$ such that $L=B\widehat L+P$ for $0\leq P< B$. 

\begin{definition}
	Let $\mathcal L$ be a language of reduced words over a set $S$. The language $\mathcal L$ has \emph{$I_P$-growth rate at least $k$} if $\mathcal L\cap R_L>ck^L$ for some $c>0$ and for all but finitely many $L\in I_P$. 
		
	Let $Q$ be a property for groups or a set of relators. Property $Q$ holds \emph{with $I_P$-overwhelming probability} if for $L\in I_P$ we have $\mathbb{P}_L(Q)\to 1$ as $L\to \infty$. 
\end{definition}

\begin{definition}
	Suppose $G=\langle S\mid  R\rangle$	is a random group with relators of length $L$. If $L\in I_0$, then any $r\in R$ can be subdivided into a word of length $\widehat L$ with letters in $\widetilde S$. We let $\widehat R$ denote the set of all such subdivided words for every $r\in R$. For $L\in I_P$, we let $\widehat R$ denote those words $q_1q_2$ such that $q_1v^{-1},vq_2\in R$ for some word $v$ of length $P$ over $S$. Then $\widehat G=\langle \widehat S \mid \widehat R\rangle$ is the \emph{associated group} to $G$ and $\widehat R$ is the set of \emph{associated relators}.
\end{definition}

\begin{lemma}
\label{lem:FCforward}
	If $\widehat G$ has property $\FC_n$	, then so does $G$. 
\end{lemma}

\begin{proof}
	Note that any $\hat r\in \widehat R$ viewed as a word over $S$ is trivial in $G$, since it either corresponds exactly to a relator or some product of relators with cancellation. Thus there exists a natural map $\widehat G\to G$. 
	
	The image $H<G$ of $\widehat G\to G$ is the subgroup generated by reduced words of length $B$ over $S$. Since $\widehat G\to H$ is surjective, $H$ has property $\FC_n$ by Lemma~\ref{lem:FCquotient}. And since $H<G$ is finite index, $G$ has property $\FC_n$ by Lemma~\ref{lem:FCvirtual}.
\end{proof}

A reduced word over $\widehat S$ is not necessarily reduced when considered as a word over $S$. However, given a c-automaton over $\widehat S$ with growth, one can retrieve a c-automaton over $S$ with a slightly reduced growth. 

\begin{definition}
	Let $\widehat\Sigma$ be a c-automaton over $\widehat S$ and $\mathcal L_{\widehat\Sigma}$ the accepted language. $\mathcal L^{\mathrm{red}}_{\widehat\Sigma}\subset \mathcal L_{\widehat\Sigma}$ is the sublanguage whose words correspond to reduced words over $S$. 	
\end{definition}

\begin{lemma}
\label{lem:reducedAutomaton}
	Let $S$ be a set of letters with $|S|=n$. Let $\widehat\Sigma$ be an automaton with $k$-growth over $\widehat S$. Then there exists an automaton $\widehat \Sigma^{\mathrm{red}}$ with $k'\defeq(k-\frac{\widehat n}{n})$-growth whose accepted language $\mathcal L_{\widehat\Sigma^{\mathrm{red}}}$ is a subset of $\mathcal L^{\mathrm{red}}_{\widehat \Sigma}$. 
	
\end{lemma}

\begin{proof}
	We begin by constructing an automaton $\widehat\Sigma'$ with the same accepted language as $\widehat\Sigma$. The vertex set $V(\widehat\Sigma')$ is the set of ordered pairs $\{(v,s) : v\in V(\widehat\Sigma)-\{v_0\}, s\in S\cup S^{-1} \}$ and a start vertex $v_0'$. For each edge $(v_1,v_2)\in E(\widehat \Sigma)$ we add following edges to $E(\widehat \Sigma')$: 
	
	\begin{enumerate}
		\item If $v_1=v_0$, then we add an edge $(v_0',(v_2,s_2))$ where the label of $(v_1,v_2)$ is an element of $\widehat S\cup \widehat S^{-1}$ which ends with $s_2$. 
		\item If $v_1\neq v_0$, then we add an edges $((v_1,s_1),(v_2,s_2))$ where $s_2$ is the last letter in the label of $(v_1,v-2)$ and $s_1\in S\cup S^{-1}$ 
	\end{enumerate}

	It is clear from construction that $\widehat\Sigma'$ and $\widehat\Sigma$ accept the same language. Additionally, the number of children is preserved, so $\Sigma$ has $k$-growth.


	To construct $\widehat\Sigma^{\mathrm{red}}$, we remove from $E(\widehat\Sigma')$ all edges of the form $((v_2,s_2),(v_3,s_3))$ whose label is a reduced word in $S$ beginning with $s_2^{-1}$. There are at most $\frac{\hat n}{n}$-many edges of this form. The resulting automaton is $\Sigma$, which has $(k-\frac{\hat n}{n})$-growth by Lemma~\ref{lem:reduceStrongGrowth}. By construction, the concatenation of any two labels on edges of the form $((v_1,s_1),(v_2,s_2))$ and $((v_2,s_2),(v_3,s_3))$ is a reduced word. Thus all words accepted by $\Sigma$ correspond to reduced words over $S$.
	
\end{proof}

\begin{corollary}
\label{cor:reducedGrowth}
	For a c-automaton $\widehat\Sigma$ with $\lambda 2\hat n$-growth, $\mathcal L^{\mathrm{red}}_{\widehat \Sigma}$ has $I_0$-growth rate at least $\sqrt[B]{\lambda-\frac{1}{2n}}(2n-1)$.
\end{corollary}

\begin{proof}
		By Lemma~\ref{lem:reducedAutomaton}, $\mathcal L^{\mathrm{red}}_{\widehat \Sigma}$ has $c\lceil 2\hat n(\lambda-\frac {1}{2n})\rceil^{\widehat L-1}\geq  c'\sqrt[B]{\big\lceil2\hat n(\lambda-\frac {1}{2n})}\big\rceil^{L}$ words of length $L\in I_0$ for some $c,c'>0$. We have
		\begin{equation*}
			\sqrt[B]{\big(\lambda - \frac {1}{2n}\big)2\hat n} > \sqrt[B]{\big(\lambda -\frac {1}{2n}\big)(2n-1)^B}=\sqrt[B]{\big(\lambda-\frac {1}{2n}\big)}(2n-1).
		\end{equation*}
\end{proof}

\begin{definition}
	Fix $1\leq P<B$ and a c-automaton $\widehat\Sigma$. We will consider $L\in I_P$. The \emph{prefix set} $\mathcal P_P\subset R_L$ is the set of reduced words whose initial length $L-P$ subword lies in $\mathcal L_{\widehat \Sigma}$. Let $u\in V(\widehat\Sigma)-\{v_0\}$ be a vertex of $\widehat\Sigma$. Let $v$ be a length $P$ reduced word over $S$ and let $s$ be the last letter of $v$. The \emph{suffix set} $\mathcal S_{u,v}$ is the set of words over $S$ of the form $vw$ such that the first letter of $w$ is not $s^{-1}$ and $w$ is an accepted word in $\widehat \Sigma_w$, the automaton $\widehat \Sigma_w$ with $T_{v_0}$ removed and $w$ set as the start vertex. 
\end{definition}

\begin{lemma}
\label{lem:suffixGrowth}
	Let $P$, $u$, and $v$ be as above. The sets $\mathcal P_P$ and $\mathcal S_{u,v}$ have $I_P$-growth rate at least $\sqrt[B]{\lambda-\frac{1}{2n}}(2n-1)$.
\end{lemma}

\begin{proof}
		Since any word in $\widehat R_{\widehat \Sigma}$ of length $L-P\in I_0$ can be completed to an element of $\mathcal P_P$, the prefix set has $I_P$-growth rate at least $\sqrt[B]{\lambda-\frac{1}{2n}}(2n-1)$ by Corollary~\ref{cor:reducedGrowth}. 
		
		$\widehat\Sigma$ has $I_0$-growth rate $\sqrt[B]{\lambda-\frac{1}{2n}}(2n-1)$ by Corollary~\ref{cor:reducedGrowth}. $\mathcal S_{u,v}$ contains all words $vw$ where $w$ is an accepted word of $\widehat \Sigma$ not beginning with $s^{-1}$. The collection of such $w$ have $I_0$-growth rate $\sqrt[B]{\lambda-\frac{1}{2n}}(2n-1)$. Consequently $\mathcal S_{u,v}$ has the desired $I_P$-growth rate.
\end{proof}

\begin{proof}[Proof of Main Theorem]
	Let $S$ be a set of generators with $n=|S|\geq 7$. Let $G=\langle S \mid R \rangle$ be a random group at density $d$ and length $L$ and let $\widehat G = \langle \widehat S \mid \widehat R \rangle$ be the associated group. Fix $B$ such that $\sqrt[B]{84}<(2n-1)^d$. Suppose $\widehat\Sigma$ is a c-automaton over $\widehat S$ with $\frac {\hat n}{6}$-growth. 
	
	Suppose that $L\in I_0$. By Corollary~\ref{cor:reducedGrowth} $\mathcal L^{\mathrm{red}}_{\widehat \Sigma}$ has $I_0$-growth rate at least $\sqrt[B]{1/84}(2n-1)$. Note $\mathcal L^{\mathrm{red}}_{\widehat \Sigma}\subset\mathcal L_{\widehat\Sigma}$ and $\sqrt[B]{1/84} (2n-1)>(2n-1)^{1-d}$ by choice of $B$. Thus by Lemma~\ref{lem:intersection} $\mathcal L_{\widehat \Sigma}\cap \widehat R$ is nonempty with $I_0$-overwhelming probability. 
	
	Suppose that $L\in I_P$ where $1\leq P<B$. By Lemma~\ref{lem:suffixGrowth}, $\mathcal P_P$ and each $\mathcal S_{u,v}$ have $I_P$-growth rate at least $\sqrt[B]{1/84} (2n-1)$. Since $\sqrt[B]{1/84} (2n-1)>(2n-1)^{1-d}$, $R$ intersects $\mathcal P_P$ and each $\mathcal S_{u,v}$ with $I_P$-overwhelming probability. Suppose that $q_1v^{-1}\in R\cap \mathcal P_P$ and suppose the path labeled by $q_1$ in $\widehat\Sigma$ ends on $u\in V(\widehat\Sigma)$. Any $w\in R\cap S_{u,v}$ is of the form $w=vq_2$. Thus word $q_1q_2$ is accepted by $\widehat\Sigma$ and belongs to $\widehat R$. Thus with $I_P$-overwhelming probability $\mathcal L_{\widehat \Sigma}\cap \widehat R$ is nonempty.
	
	We have shown that with overwhelming probability $\mathcal L_{\widehat \Sigma}\cap \widehat R$ is nonempty for any c-automaton $\widehat\Sigma$ with $\frac {\hat n}{6}$-growth. We claim $\widehat G$ has $\FC_2$. By Theorem~\ref{thm:main1} there exists a finite collection of c-automata $\widehat \Sigma^1,\dots, \widehat \Sigma^k$ with $\frac{\hat n}{6}$-growth such that for any $\mathrm{CAT}(0)$ square complex $\X$ and action $F_{\widehat S}\to \Aut(\X)$ without global fixed point, every $w\in \mathcal L_{\widehat \Sigma^j}$ acts nontrivially on $\X$ for some $j\in\{1,\ldots, k\}$. Suppose for contradiction that some $F_{\widehat S}\to \Aut(\X)$ induces an action of $\widehat G$ on a $\mathrm{CAT}(0)$ square complex $\widetilde X$ without global fixed point. We get a contradiction since $\mathcal L_{\widehat \Sigma}\cap \widehat R$ is nonempty w.o.p. By Lemma~\ref{lem:FCforward}, $G$ has $\FC_2$.
\end{proof}

\bibliographystyle{amsalpha}
\bibliography{bibliography}

\end{document}